\documentstyle[twoside,psfig]{article}

\oddsidemargin=\evensidemargin
\catcode`\@=11
\long\def\@makefntext#1{
\protect\noindent \hbox to 3.2pt {\hskip-.9pt  
$^{{\eightrm\@thefnmark}}$\hfil}#1\hfill}		

\def\ps@myheadings{\let\@mkboth\@gobbletwo		
\def\@oddhead{\hbox{}
\rightmark\hfil\eightrm\thepage}   
\def\@oddfoot{}\def\@evenhead{\eightrm\thepage\hfil
\leftmark\hbox{}}\def\@evenfoot{}
\def\sectionmark##1{}\def\subsectionmark##1{}}

\catcode`@=11		      		     
\def\ps@plain{\let\@mkboth\@gobbletwo
     \def\@oddhead{}\def\@oddfoot{\eightrm\hfil\thepage
     \hfil}\def\@evenhead{}\let\@evenfoot\@oddfoot}

\newcounter{sectionc}\newcounter{subsectionc}\newcounter{subsubsectionc}
\renewcommand{\section}[1] {\vspace{12pt}\addtocounter{sectionc}{1} 
\setcounter{subsectionc}{0}\setcounter{subsubsectionc}{0}\noindent 
	{\tenbf\thesectionc. #1}\par\vspace{5pt}}
\renewcommand{\subsection}[1] {\vspace{12pt}\addtocounter{subsectionc}{1} 
	\setcounter{subsubsectionc}{0}\noindent 
	{\bf\thesectionc.\thesubsectionc. 
	{\kern1pt \bfit #1}}\par\vspace{5pt}}
\renewcommand{\subsubsection}[1] {\vspace{12pt}
	\addtocounter{subsubsectionc}{1}
	\noindent
	{\tenrm\thesectionc.\thesubsectionc.\thesubsubsectionc.	{\kern1pt 
	\it #1}}\par\vspace{5pt}}

\newcounter{appendixc}
\newcounter{subappendixc}[appendixc]
\newcounter{subsubappendixc}[subappendixc]

\renewcommand{\appendix}[1] {\vspace{12pt}	
	\refstepcounter{appendixc}		
	\setcounter{figure}{0}
	\setcounter{table}{0}
	\setcounter{lemma}{0}
	\setcounter{theorem}{0}
	\setcounter{corollary}{0}
	\setcounter{definition}{0}
	\setcounter{equation}{0}
	\renewcommand{\thefigure}{\Alph{appendixc}.\arabic{figure}}
	\renewcommand{\thetable}{\Alph{appendixc}.\arabic{table}}
	\renewcommand{\theappendixc}{\Alph{appendixc}}
	\renewcommand{\thelemma}{\Alph{appendixc}.\arabic{lemma}}
	\renewcommand{\thetheorem}{\Alph{appendixc}.\arabic{theorem}}
	\renewcommand{\thedefinition}{\Alph{appendixc}.\arabic{definition}}
	\renewcommand{\thecorollary}{\Alph{appendixc}.\arabic{corollary}}
	\renewcommand{\theequation}{\Alph{appendixc}.\arabic{equation}}
	\noindent{\tenbf Appendix \theappendixc #1}\par\vspace{5pt}}

\newcommand{\textlineskip}{\baselineskip=13pt}
\newcommand{\smalllineskip}{\baselineskip=10pt}

\def\abstracts#1#2#3#4{{
	\centering{\begin{minipage}{4.5in}\footnotesize\baselineskip=10pt
	\centerline{ABSTRACT} 
	\parindent=15pt #1\par 
	\parindent=15pt #2\par
	\parindent=15pt #3\par
	\parindent=15pt #4\par
	\end{minipage}}\par}} 

\def\keywords#1{{ 
	\centering{\begin{minipage}{4.5in}\footnotesize\baselineskip=10pt
	{\footnotesize\it Keywords}\/: #1
	\end{minipage}}\par}}

\renewenvironment{thebibliography}[1]
	{\frenchspacing
	 \ninerm\baselineskip=11pt
	 \begin{list}{[\arabic{enumi}]}
	{\usecounter{enumi}\setlength{\parsep}{0pt}
	 \setlength{\leftmargin 13.7pt}{\rightmargin 0pt} 
	 \setlength{\itemsep}{0pt} \settowidth
	{\labelwidth}{[#1]}\sloppy}}{\end{list}}

\newcounter{itemlistc}
\newcounter{romanlistc}
\newcounter{alphlistc}
\newcounter{arabiclistc}

\newcommand{\fcaption}[1]{
        \refstepcounter{figure}
        \setbox\@tempboxa = \hbox{\footnotesize Fig.~\thefigure. #1}
        \ifdim \wd\@tempboxa > 5in
           {\begin{center}
        \parbox{5in}{\footnotesize\smalllineskip Fig.~\thefigure. #1}
            \end{center}}
        \else
             {\begin{center}
             {\footnotesize Fig.~\thefigure. #1}
              \end{center}}
        \fi}

\newcommand{\tcaption}[1]{
        \refstepcounter{table}
        \setbox\@tempboxa = \hbox{\footnotesize Table~\thetable. #1}
        \ifdim \wd\@tempboxa > 5in
           {\begin{center}
        \parbox{5in}{\footnotesize\smalllineskip Table~\thetable. #1}
            \end{center}}
        \else
             {\begin{center}
             {\footnotesize Table~\thetable. #1}
              \end{center}}
        \fi}

\def\pmb#1{\setbox0=\hbox{#1}
	\kern-.025em\copy0\kern-\wd0
	\kern.05em\copy0\kern-\wd0
	\kern-.025em\raise.0433em\box0}

\def\fnt#1#2{\footnotetext{\kern-.3em
	{$^{\mbox{\scriptsize #1}}$}{#2}}}

\def\fpage#1{\begingroup
\voffset=.3in
\thispagestyle{empty}\begin{table}[b]\centerline{\footnotesize #1}
	\end{table}\endgroup}

\def\runninghead#1#2{\pagestyle{myheadings}
\markboth{{\protect\footnotesize\it{\quad #1}}\hfill}
{\hfill{\protect\footnotesize\it{#2\quad}}}}

\font\tenrm=cmr10
 
\font\tenbf=cmbx10
\font\bfit=cmbxti10 at 10pt
\font\ninerm=cmr9

\font\eightrm=cmr8

\newtheorem{theorem}{Theorem}   

\newtheorem{corollary}{Corollary}
\def\@begintheorem#1#2{\trivlist	
	\item[\hskip\labelsep{\bf #1\ #2.}]} 
\def\@opargbegintheorem#1#2#3{\trivlist
	\item[\hskip\labelsep{\bf #1\ #2\ (#3).}]}

    	{\setcounter{itemlistc}{0}		
	 \begin{list}{$\bullet$}		
	{\usecounter{itemlistc}			
	 \leftmargin10pt	       
	 \setlength{\parsep}{0pt}
	 \setlength{\itemsep}{0pt}     
	}}{\end{list}}

	{\setcounter{romanlistc}{0}		
	 \begin{list}{$($\roman{romanlistc}$)$}	
	{\usecounter{romanlistc}		
	 \leftmargin18pt 
	 \setlength{\parsep}{0pt}
	 \setlength{\itemsep}{0pt}	
	 \settowidth{\labelwidth}{#1}                          
	}}{\end{list}}

	{\setcounter{enumii}{0}			
	 \begin{list}{$($\alph{enumii}$)$}	
	{\usecounter{enumii}			
	 \leftmargin18pt		
	 \setlength{\parsep}{0pt}
	 \setlength{\itemsep}{0pt}	
	 \settowidth{\labelwidth}{#1}                          
	}}{\end{list}}

\def\qed{\hbox{${\vcenter{\vbox{			
   \hrule height 0.4pt\hbox{\vrule width 0.4pt height 6pt
   \kern5pt\vrule width 0.4pt}\hrule height 0.4pt}}}$}}

\def\theequation{\thesectionc.\arabic{equation}}  
\newtheorem{example}{Example}

\def \Q {{\bf Q}}
\def \Z {{\bf Z}}

\def\AA {\ifmmode{{\cal A}}\else{${\cal A}$}\fi}
\def\TT {\ifmmode{{\cal T}}\else{${\cal T}$}\fi}
\def\LL {\ifmmode{{\cal L}}\else{${\cal L}$}\fi}
\def\subl{\hookrightarrow}
\def\divs{\ \vrule height12pt  depth6pt\ }
\def\Ma{\ifmmode{M(\alpha)}\else{$M(\alpha)$}\fi}
\def\Mb{\ifmmode{M(\beta)}\else{$M(\beta)$}\fi}

\begin{document}
\setlength{\textheight}{7.7truein}  

\runninghead{Embedding tangles in links}
{Embedding tangles in links}

\normalsize\textlineskip
\thispagestyle{empty}
\setcounter{page}{1}


\vspace*{0.88truein}
\fpage{1}
\centerline{\bf EMBEDDING TANGLES IN LINKS}
\baselineskip=13pt
\vspace*{0.37truein}
\centerline{\footnotesize DANIEL RUBERMAN}
\baselineskip=12pt
\centerline{\footnotesize\it Department of Mathematics}
\baselineskip=10pt
\centerline{\footnotesize\it MS 050}
\centerline{\footnotesize\it Brandeis University}
\centerline{\footnotesize\it Waltham, MA 02254}

\vspace*{10pt}

%
\abstracts{Extending and reproving a recent result of D.~Krebes, we give
obstructions to the embedding of a tangle in a link.}{}{}{}

\vspace*{10pt}


\vspace*{1pt}\textlineskip	
\section{Introduction}	
\vspace*{-0.5pt}

A recent paper~\cite{krebes:tangles} of David Krebes poses the following
interesting question: given a tangle \TT, and a link \LL, when can \TT\ sit
inside \LL?  A {\sl tangle} is a $1$-manifold with 4 boundary components,
properly embedded in a 3-ball; this is somewhat more general than the usual
definition.  An embedding of \TT\ in \LL\  is determined by a ball in $S^3$,
whose intersection with \LL\  is the given tangle; we will indicate
an embedding
by $\TT \subl \LL$.  Krebes gives the following simple criterion for when $\TT$
embeds in \LL.  Complete the tangle to a link in either of the two obvious
ways, giving rise to the {\sl numerator}
$n(\TT)$ and {\sl denominator} $d(\TT)$.  For any oriented
link $\LL$, its determinant
$\det(\LL)$ is defined to be $\det(V + V')$, where $V$ is a Seifert
matrix for $\LL$.   Krebes shows
\begin{theorem}\label{Kther}\sl
If $\TT \subl \LL$, then
\begin{equation}
\gcd(\det(n(\TT)),\det(d(\TT)))\divs \det(\LL)
\end{equation}
\rm
\end{theorem}

The proof of Theorem~\ref{Kther} given in~\cite{krebes:tangles} is essentially
combinatorial, and uses an interpretation of the determinant in terms of link
diagrams and the Kauffman bracket.  On the other hand, the determinant has
a homological interpretation; it is essentially the order of the
homology of the
$2$-fold branched cover $S^3$, branched along the link.   In this paper, we
will prove a simple fact about the homology of certain $3$-manifolds, which
readily implies Theorem~\ref{Kther}.   In essence, we replace the
divisibility condition above with the conclusion that the homology of the
$2$-fold branched cover of the link must contain a subgroup of a certain size.
Our approach has the advantage that it gives some stronger results on the
embedding problem, which do not seem approachable via the combinatorial route.
A slightly different argument yields a similar conclusion about other branched
coverings as well. Some examples, and further remarks on embeddings, are
given in the final section.

To state the result, let us use the notation $|M|$ for the order of the first
homology group of $M$, where by convention the order is defined to be $0$ if
the homology is infinite.  Also, if $\partial M$ is a torus $T^2$, and $\alpha
\subset T$ is a simple closed curve which does not bound a disc, then let
$M(\alpha)$ denote the result of Dehn filling with slope $\alpha$.  (In other
words, glue $S^1 \times D^2$ to $M$ so that $\partial D^2$ is glued to
$\alpha$.)
\begin{theorem}\label{Rther}\sl
Suppose that $M$ is an orientable $3$-manifold, and that $\alpha,\beta$ are
simple closed curves on $T = \partial M$ which generate all of $H_1(T)$.
Suppose that $M \subset N$, where $N$ is a closed, orientable $3$-manifold.
Then
$$
\gcd(|\Ma|,|\Mb|) \divs |N|
$$
\rm
\end{theorem}
It is worth remarking that as a consequence of Theorem~\ref{Rther}, the
quantity 
$$f(M) = \gcd(|\Ma|,|\Mb|)$$
is independent of the choice of pair $\alpha,
\beta$, and hence defines an invariant of $M$.  For, given another such pair,
say
$\alpha',\beta'$, the theorem says that $\gcd(|\Ma|,|\Mb|)$ divides both
$|M(\alpha')|$ and $|M(\beta')|$, and so  $\gcd(|\Ma|,|\Mb|) \divs
\gcd(|M(\alpha')|,|M(\beta')|)$.  The remark follows by reversing the roles of
$\alpha,\beta$ and $\alpha',\beta'$.

Some further results on the embedding problem, using invariant derived from the
Kauffman bracket, can be found in the
recent preprint~\cite{krebes-silver-williams:tangles}.

\vspace*{1pt}\textlineskip	
\section{Proof of Theorem~\protect{\ref{Rther}}}
\vspace*{-0.5pt}

Before beginning the proof of the theorem, we remark that unless $H_1(N)$ is
torsion, then the theorem has no content.  So we can assume that $N$
is a rational homology sphere for the remainder of this section.  Writing
$$
N = M \cup_T X
$$
it follows from a standard Poincar\'e duality argument that both $M$ and $X$
have the rational homology of a circle.  Hence we can write (non-canonically)
$$
H_1(M) = \Z \oplus \Z/q_1 \oplus \dots \Z/q_s.
$$
In particular, the torsion subgroup $T_1(M) \subset H_1(M)$ has order
$q_1\cdots q_s$.

Under the map $j_*:H_1(T) \to H_1(M) $, the classes $\alpha $ and $\beta$ go
to
$$
j_*(\alpha) = (a,a_1,\dots,a_s)\quad {\rm and} \quad j_*(\alpha) =
(b,b_1,\dots,b_s)
$$
respectively.  (The coefficients are with respect to generators of the summands
of $H_1(M)$ in the splitting given above.)\\[1ex]
{\bf Claim:} $|\Ma| = a\,|T_1(M)|.$\\[1ex]
To see this, note that $H_1(\Ma) \cong H_1(M)/<\alpha>$. The given splitting of
$H_1(M)$ corresponds to a presentation of that group by the $s \times (s+1)$
matrix
\begin{equation}
\pmatrix{
0& q_1 & 0 & \dots& 0\cr
0&0 & q_2  & \dots& 0\cr
\dots & \dots & \dots & \dots& \dots\cr
0&0 & 0  & \dots& q_s\cr
}
\end{equation}
Killing $\alpha$ adds an additional row, to get the following presentation
matrix for
$H_1(\Ma)$:
\begin{equation}
\pmatrix{
a& a_1 & a_2 & \dots& a_s\cr
0& q_1 & 0 & \dots& 0\cr
0&0 & q_2  & \dots& 0\cr
\dots & \dots & \dots & \dots& \dots\cr
0&0 & 0  & \dots& q_s\cr
}
\end{equation}
which has determinant $aq_1 \cdots q_s =  a\,|T_1(M)|.$  But the order of
$H_1(M)$ is the same as the determinant of any (square) presentation matrix for
it.

By the same argument for homology of $\beta $ we get that
\begin{equation}
\gcd(|\Ma|,|\Mb|) = \gcd (a\, |T_1(M)|, b\, |T_1(M)|) = \gcd(a,b) |T_1(M)|
\end{equation}

Now we turn to the situation at hand, and consider the homology of $N$, which we
place into the long exact sequence of the pair $(N,M)$.
$$
\begin{array}{ccccccccccc}
0&\longrightarrow & H_2(N,M) & \stackrel{\partial}{\longrightarrow} & H_1(M) 
& \longrightarrow& H_1(N)
&\longrightarrow H_1(N,M) 
&\longrightarrow &0\\
&&\Big\uparrow\vcenter{%
\rlap{$\scriptstyle \cong$}}
&&
\Big\uparrow\vcenter{%
\rlap{$\scriptstyle j_*$}}
&&&&&\\
&&H_2(X,T)& \stackrel{\partial}{\longrightarrow} &H_1(T)&&&&&
\end{array}
$$

By exactness, $H_1(M)/\partial(H_2(N,M))$ injects into $H_1(N)$, and so the
theorem will follow if we can show that $\gcd(a,b) |T_1(M)|$ divides the order
of  $H_1(M)/\partial(H_2(N,M))$.   Again by duality,
$H_2(X,T) \cong H^1(X)
\cong
\Z$ is generated by a relative 2-cycle $C$ with boundary $\gamma$ lying in $T$.
By the commutativity of the diagram (the isomorphism is an excision) it follows
that $H_1(M)/\partial(H_2(N,M))$ is just $H_1(M)/<\gamma>$.
As before, the image of $\gamma$ in $H_1(M)$ may be written
$ (c,c_1,\dots,c_s)$, so that $|H_1(M)/<\gamma>| = c\, |T_1(M)|$.   Since
$\alpha$ and $\beta$ generate $H_1(T) $, we can write  $\gamma = m\alpha
+n \beta$ which implies that $c = ma +nb$.  This means that
$\gcd(a,b) \divs c$,
or in other words
$$
\gcd(a,b)|T_1(M)| \divs c\, |T_1(M)|.
$$
Since $c\, |T_1(M)|\divs |H_1(N)|$ the theorem follows.\hfill{}\qed

\vspace*{1pt}\textlineskip	
\section{Complements and Examples}
\vspace*{-0.5pt}

First, let us observe that Theorem~\ref{Rther} implies Theorem~\ref{Kther}.
The basic point is that the $2$-fold cover of the ball, branched along a
trivial
tangle, is $S^1 \times D^2$.  The different ways of completing a tangle
\TT\  to
a link give rise to different Dehn fillings of $M$ = the $2$-fold cover of
the ball, branched along \TT.  It is easy to see that the
meridians of the solid tori corresponding to the numerator $n(\TT)$ and
denominator $d(\TT)$ have intersection number $\pm1$ in $T=\partial M$ and
thus generate $H_1(T)$. Now if
$\TT \subl \LL$ as in Theorem~\ref{Kther}, there is an embedding
$M \subset N$, where $N$ is the $2$-fold cover of the $3$-sphere branched
along \LL.  The $2$-fold cover of the ball, branched along a
trivial tangle, is a solid torus.  Hence the hypotheses of Theorem~\ref{Rther}
are satisfied, and we get that
$$
\gcd(| n(\TT)|,| d(\TT)|) \divs |\det(\LL)|
$$

Second, the proof gives a somewhat stronger condition, involving the homology
groups themselves, rather than their orders.  We give the statement in terms of
the homology of $3$-manifolds, with the understanding that passing to the
2-fold cover of a link gives rise to a restriction on embeddings of tangles in
links.
\begin{corollary}\label{summands}\sl
Suppose that $H_1(N)$ is torsion, and that $M \subset N$.  Then there is an
injection of $\Z/c \oplus T_1(M)$ into $H_1(N)$, where $c$ has the same
significance as in the proof of Theorem~\ref{Rther}.
\rm
\end{corollary}

In applying these results to specific tangles, the most useful part of the
conclusion is the fact 
that the inclusion map of $M$ into
$N$ induces an injection on $T_1(M)$.  This is true in a more general setting,
by an argument which is perhaps more conceptual than the calculation proving
Theorem~\ref{Rther}.
\begin{theorem}\label{torsion}\sl
Suppose $M$ is an orientable $3$-manifold with connected boundary, and $i: M
\subset N$, where $N$ is an orientable $3$-manifold with $H_1(N)$ torsion.
Then
the inclusion map $i_*$ induces an injection of $T_1(M) $ into $H_1(N)$.
\rm
\end{theorem}
\noindent
{\it Proof of Theorem~\ref{torsion}:}
The linking pairing $\lambda: T_1(M) \times T_1(M,\partial M) \to \Q/\Z$ is
non-degenerate, by Poincar\'e duality. So if $x \in T_1(M)$ is non-zero, there
is an element $y \in  T_1(M,\partial M)$ with $\lambda(x,y)$ non-zero in
$\Q/\Z$; represent each of these by (absolute or relative) cycles with the
same names.  Since
$\partial M$ is connected, $H_1(M)
\to H_1(M,\partial M)$ is surjective, so $y$ lifts to $\bar{y} \in H_1(M)$.
There is no good reason that $\bar{y}$ represents a torsion class in $M$, but
by hypothesis it is a torsion class in $N$.  Now $\lambda(x,i_*\bar{y})$ (as
calculated in $N$) may be calculated as the intersection number
of $\bar{y}$ with a $2$-chain bounding $n\cdot x$, and $C$ can be chosen to lie
in $M$ since $x$ is a torsion class in $M$.  Hence
$\lambda(x,i_*\bar{y})= \lambda(x,y) \neq 0$ in $\Q/\Z$, and therefore $x$ is
nontrivial in $H_1(N)$.\hfill{}\qed\\[1ex]
{\bf Remark:}
It is not possible to deduce
Theorem~\ref{Rther} (and hence Theorem~\ref{Kther}) from
Theorem~\ref{torsion}.  To do so would amount to proving that (in the notation
of Theorem~\ref{Rther})
$\gcd(|\Ma|,|\Mb|) = |T_1(M)|$.  However, this is not the case, as the
following example indicates; the example has $T_1(M) = \Z/3$ but
$\gcd(|\Ma|,|\Mb|) =9.$

Consider an oriented solid torus $M_0$  with a basis $\alpha,\beta$ for $H_1(T)$
chosen so that $\beta$ generates $H_1(M_0)$ and $\alpha$ bounds a disk.  Let
$K \subset M_0$ be an oriented knot, representing
$3$ times $\beta$ in $H_1(M_0)$.  The meridian $\mu$ of $K$ is determined by
the orientation, and we choose a longitude $\lambda$ by requiring that
$\lambda$ be homologous to $3 \beta$ in $H_1(M_0 - K)$. (If $M_0$ were embedded
in $S^3$ in a standard way, so that $\beta$ bounds a disk in $S^3 - {\rm
int}(M_0)$, then $\lambda$ would be the longitude of $K$ in $S^3$.)  Now let
$M$ be the result of Dehn surgery on $M_0$, with coefficient $9/1$.  In other
words, remove a neighborhood of $K$, and glue in a solid torus killing $9\mu +
\lambda$.  The homology of $M_0 - \nu(K)$ is generated by $\mu$ and $\beta$,
and the surgery kills $9\mu + 3\beta$, so $H_1(M)= \Z + \Z/3$.  

On the other hand, the homology of $\Ma$ is gotten by killing $\alpha$, which
is homologous to $3\mu$ and so $H_1(\Ma)$ is presented by the matrix
\begin{equation}
\pmatrix{
9& 3\cr
3&0\cr
}
\end{equation}
yielding $H_1(\Ma) = \Z/3 \oplus \Z/3$.  Likewise, $H_1(\Mb)$ is presented by
the matrix
\begin{equation}
\pmatrix{
9& 3\cr
0&1\cr
}
\end{equation}
yielding $H_1(\Mb) = \Z/9$.

This same example shows that the hypothesis in Theorem~\ref{torsion}, that
$H_1(N)$ be torsion, is necessary.  For if $N$ is obtained by filling $M$ with
slope $\alpha + \beta$, then $H_1(N) \cong \Z$, and so $T_1(M)$ doesn't
inject.\\[1ex]

Theorem~\ref{torsion} may be applied to branched covers of tangles, of degrees
other than $2$.  In applying this remark, one must take care, because for
$k> 2$
the $k$-fold covers of the ball, branched along \TT\ are not uniquely specified
by the branch locus.  These different covers may well have differing homology
groups, so in practice, one might have to calculate the homology groups for all
of the different possibilities.
\begin{corollary}\label{kfold}\sl
Suppose that $\TT \subl \LL$, and that $N$ is a k-fold cover of $S^3$
branched along \LL.  Let $M$ be the induced cover of $B^3$, branched along
\TT.   If $N$ is a rational homology sphere, then
$T_1(M)$ is a subgroup of $T_1(N)$.\rm
\end{corollary}

Most of the results so far have concerned only the torsion part of the
homology, but there are some things which can be said about the torsion-free
part of the homology.   One simple result is the following.
\begin{theorem}\label{rational}\sl
Suppose that $M$ has boundary of genus $g$, and that $i: M \subset N$.  Then
$\dim(H_1(N;\Q)) \ge \dim(H_1(M;\Q)) -g$.
\rm
\end{theorem}

The proof is straightforward; if the quantity $\dim(H_1(M;\Q)) -g$ is greater
than zero, then there is a subspace of that dimension in $H_1(M;\Q)$ which
pairs non-trivially with a subspace of $H_2(M;\Q)$ of the same dimension.
Hence both of these inject into the homology of $N$.
\begin{example}
It is not hard to give examples where the homological approach gives more
information than can be deduced from the determinants alone.  The simplest
I can
think of is the following.  Let $\TT$ be the tangle
$$
(T_3)^* + (T_3)^* + (T_{-3})^*
$$
pictured below in Figure 1.\\
\begin{figure}[htbp] 
\vspace*{1pt}
\centerline{\psfig{file=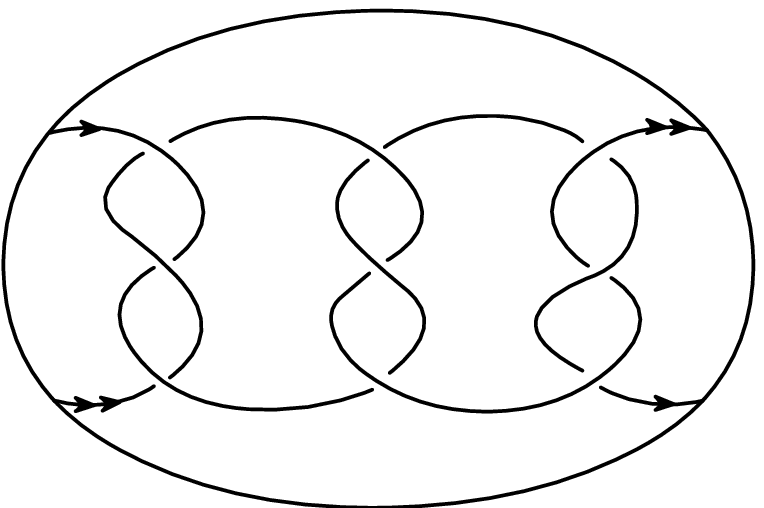,height=1.9in}} 
\vspace*{2pt}
\centerline{\bf Figure 1}
\end{figure}
According to~\cite{krebes:tangles}, the numerator
and denominator of this tangle may be calculated as the numerator and
denominator of the fraction obtained by the grade school addition of fractions,
without canceling common factors:
$$
\frac{-1}{3} + \frac{-1}{3} + \frac{1}{3} = \frac{-6}{9} + \frac{1}{3} =
\frac{-18 + 9}{27} = \frac{-9}{27}
$$
and so $\gcd(|n(\TT)|,|d(\TT)|) = 9$.  This would allow, in principle, that
$\TT$ might embed in a $2$-bridge knot (or link) corresponding to the fraction
$9p/q$, for the determinant of such a knot is $9p$.  But we calculate below
that
for
$M = $ the $2$-fold cover of $B^3$ branched along $\TT$, we have $H_1(M) = \Z
\oplus \Z/3 \oplus \Z/3 $, whereas the $2$-fold cover of a $2$-bridge knot or
link has cyclic homology.
\end{example}

The calculation of $H_1(M)$ proceeds as in the usual calculation of $2$-fold
covers branched along knots, as described in \S6 of~\cite{rolfsen:knots}.  If
surgery is done on the middle crossings in each of the three $\pm 1/3$ tangles
which make up \TT, then \TT\ becomes trivial.  Straightening it out (for the
purpose of drawing the branched cover, it's legal to move the endpoints
around) gives the surgery picture in Figure 2.
\begin{figure}[htbp] 
\vspace*{1pt}
\centerline{\psfig{file=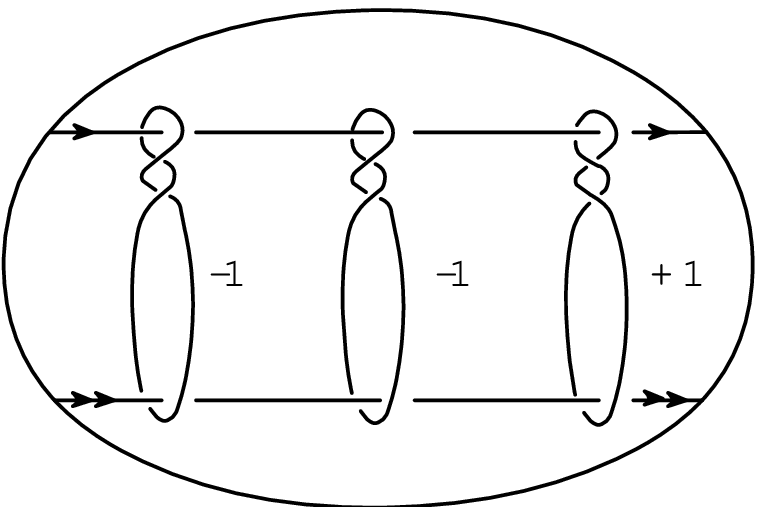,height=1.9in}} 
\vspace*{2pt}
\centerline{\bf Figure 2}
\end{figure}

Passing to the double branched cover give a picture of $M$ as surgery on a
$6$-component link in $S^1 \times D^2$, whose homology may be readily
calculated to give the result quoted above.

\vspace*{1pt}\textlineskip
\section{An obstruction to embedding in a trivial link}
\vspace*{-0.5pt}

An obstruction to embedding tangles in the trivial link, of a somewhat
different sort, may be derived from the invariants $I^n(\TT)$ defined
in~\cite{cochran-ruberman:tangles}.  To explain this, we recall that for a
$2$-component oriented link $\LL = (L_x,L_y)$ of linking number $\lambda =0$,
Cochran~\cite{cochran:gilc} defined `higher linking numbers' $\beta^n_x(\LL)$
and
$\beta^n_y(\LL)$.  For $n=1$ these are both equal to the Sato-Levine
invariant~\cite{sato} while for $n>1 $ they depend on the ordering of the
components.  For a tangle $\TT$ it is possible to choose a tangle sum with a
trivial tangle (a {\sl closure} of \TT\ in the terminology
of~\cite{cochran-ruberman:tangles}) to get a link $\LL$ with $\lambda= \beta^1
= 0$.  There is no canonical choice for \LL, but we showed that
$$
|I^n(\TT)| = |\beta^n_x(\LL)- \beta^n_y(\LL)|
$$
independent of $\LL$ (and of the order of the components because of the
absolute value signs).
\begin{theorem}\label{cr}
\sl
Suppose that $\TT$ is a tangle with no
loops such that $\TT \subl {\cal J} = $ the trivial $2$-component link.  Then
$I^n(\TT) = 0$ for all $n\ge 2$.
\rm
\end{theorem}
\noindent
{\it Proof:} Consider first a $2$-string tangle ${\cal K}$ with no
loops, and
the 4-punctured sphere $\Sigma$ in the boundary of its exterior.  Then $\Sigma$
is compressible if and only if ${\cal K}$ is {\sl split}, where
${\cal K}$
is split if the ball can be split into two sub-balls by a properly embedded
disk with the two strings of ${\cal K}$ lying in different sub-balls.   This
follows directly from the loop theorem and Dehn's lemma.  The structure of a
split tangle is very simple: it is a trivial tangle possibly with knots tied in
each string.  The tangle is trivial if and only if there are no knots.  In
particular, a split tangle has $I^n = 0$ for all $n$, because it has a
completion which is a split link, which in turn has all of its $\beta^n = 0$.

Now consider a tangle \TT\ with $I^n \neq 0$, and suppose that $\TT \subl
{\cal J}$, or in other words that ${\cal J}$ splits as a sum $\TT +
\TT'$.  From the preceding paragraph $I^n(\TT)
\neq 0$ implies that the surface $\Sigma$ is incompressible in the exterior of
\TT.  But then $\Sigma$ must be compressible in the exterior of $\TT'$. For if
it weren't then $\Sigma$ would be an incompressible $4$-punctured sphere in the
exterior of the unlink.  Applying the preceding paragraph once more, it
follows
that $\TT'$ is split.  If there is a knot in one of the strands of $\TT'$, then
that would give a knot in the corresponding component of the unlink, which
cannot be.  It follows from all of this that $\TT'$ must in fact be a trivial
tangle, so that the the unlink $\TT + \TT'$ may be used to calculate $I^n(\TT)$
and show that it is zero.  This contradicts our assumption that $I^n(\TT)
\neq 0$.\hfill{}\qed\\[1ex]

The tangles cited in~\cite{cochran-ruberman:tangles} give rise to examples of
tangles which cannot be embedded in a trivial link.

\vspace*{1pt}\textlineskip	
\section{Generalizations of Tangles}
\vspace*{-0.5pt}

We close with a few remarks on some of the questions raised at the end
of~\cite{krebes:tangles}.  One question concerned the existence of a family of
completions of any tangle with $2t$ strands, which would play the role of the
numerator and denominator of a $2$-string tangle.  Following the proof of
Theorems~\ref{Rther} and~\ref{torsion} we see how to construct such a family.
Note that the $2$-fold branched cover of a trivial $2t$-tangle is a handlebody
of genus $2t-1$.  By twisting the strings around, one can vary the attaching of
this complementary handlebody so as to kill of the homology of $\partial
M$ in various fillings.   The $\gcd$ of the orders of the homology of all these
fillings then would have to divide the determinant of any link in which the
tangle was embedded.  

Another generalization of a tangle pointed out in \S14 of~\cite{krebes:tangles}
is an embedded arc (or more generally a $1$-manifold) in a more complicated
submanifold of $S^3$, such as a solid torus.  In particular,
Krebes asks whether a particular arc ${\cal A}$ in a $S^1 \times D^2$ (Figure 17
of~\cite{krebes:tangles})  can sit inside an unknot.
Our approach, especially Corollary~\ref{kfold}
gives (in principle) an obstruction, if there were torsion in the
homology of some cover of $S^1 \times D^2$ branched along ${\cal A}$.
Unfortunately, it seems that the homology of all of the cylic covers of the
solid torus, branched along this arc, is torsion-free.  Hence we cannot apply
Theorem~\ref{torsion} to deduce anything about embeddings of this pair in a
link.  Likewise, it does not seem possible to use Theorem~\ref{rational},
because the rational homology of each of the cyclic branched covers is the same
as for a trivial arc.  

\section{Acknowledgements}
The author was partially supported by NSF grant DMS 9917802.

\end{document}